\documentclass[11pt]{article}
\usepackage{epsfig}

\def\be{\begin{equation}}
\def\ee{\end{equation}}
\def\bea{\begin{eqnarray}}
\def\eea{\end{eqnarray}}
\def\bes{\begin{eqnarray*}}
\def\ees{\end{eqnarray*}}

\def\nn{\nonumber}
\def\<{\langle}
\def\>{\rangle}
\def\lb{\label}
\def\bs{\setminus}

\def\R{{\bf R}}
\def\C{{\bf C}}
\def\Z{{\bf Z}}
\def\N{{\bf N}}
\def\U{{\bf U}}
\def\Q{{\bf Q}}
\def\T{{\bf T}}

\def\aa{{\alpha}}

\def\ga{{\gamma}}

\def\th{{\theta}}
\def\om{{\omega}}
\def\Om{{\Omega}}
\def\ep{{\epsilon}}
\def\lm{{\lambda}}

\def\sg{{\sigma}}

\def\Sg{{\Sigma}}
\def\vf{{\varphi}}

\def\H{{\cal H}}
\def\T{{\cal T}}

\def\P{{\cal P}}

\def\Nn{{\cal N}}

\def\Sp{{\rm Sp}}

\def\dm{{\rm \diamond}}

\def\ol#1{\overline{#1}}  

\def\hb{\vrule height0.18cm width0.14cm $\,$}

\def\ol#1{\overline{#1}}  

\title{Closed trajectories on symmetric convex \\Hamiltonian energy surfaces}

\author{Wei Wang\thanks{Partially supported by National Natural
Science Foundation of China No.10801002, China Postdoctoral Science
Foundation No.20070420264 and LMAM in Peking University.
E-mail: alexanderweiwang@yahoo.com.cn, wangwei@math.pku.edu.cn  }\\
School of Mathematical Science \\ Peking University, Beijing 100871 \\
PEOPLES REPUBLIC OF CHINA \\ }

\date{July 4th, 2009}
\begin{document}

\maketitle

\begin{abstract}
{\it In this article, let $\Sigma\subset\R^{2n}$ be a compact convex
Hamiltonian energy surface which is symmetric with respect to the
origin. where $n\ge 2$. We prove that there exist at least two
geometrically distinct symmetric closed trajectories of the Reeb
vector field on $\Sg$. }
\end{abstract}



\renewcommand{\theequation}{\thesection.\arabic{equation}}
\renewcommand{\thefigure}{\thesection.\arabic{figure}}

\setcounter{equation}{0}
\section{Introduction and main results}

In this article, let $\Sigma$ be a fixed $C^3$ compact convex hypersurface
in $\R^{2n}$, i.e., $\Sigma$ is the boundary of a compact and strictly
convex region $U$ in $\R^{2n}$. We denote the set of all such hypersurfaces
by $\H(2n)$. Without loss of generality, we suppose $U$ contains the origin.
We denote the set of all compact convex hypersurfaces which are
symmetric with respect to the origin by $\mathcal{SH}(2n)$,
i.e., $\Sg=-\Sg$ for $\Sg\in\mathcal{SH}(2n)$.
We consider closed characteristics $(\tau,y)$ on $\Sigma$, which are
solutions of the following problem
\be \left\{\matrix{\dot{y}=JN_{\Sigma}(y), \cr
               y(\tau)=y(0), \cr }\right. \lb{1.1}\ee
where $J=\left(\matrix{0 &-I_n\cr
                I_n  & 0\cr}\right)$,
$I_n$ is the identity matrix in $\R^n$, $\tau>0$ and $N_\Sigma(y)$ is
the outward normal vector of $\Sigma$ at $y$ normalized by the
condition $N_{\Sigma}(y)\cdot y=1$. Here $a\cdot b$ denotes the
standard inner product of $a, b\in\R^{2n}$. A closed characteristic
$(\tau,\, y)$ is {\it prime} if $\tau$ is the minimal period of $y$.
Two closed characteristics $(\tau,\, y)$ and $(\sigma, z)$ are {\it
geometrically distinct}  if $y(\R)\not= z(\R)$. We denote by
$\T(\Sg)$ the set of geometrically distinct closed
characteristics $(\tau,\, y)$ on $\Sg$.
A closed characteristic $(\tau, y)$ on $\Sg\in \mathcal{SH}(2n)$ is
{\it symmetric} if $\{y(\R)\}=\{-y(\R)\}$,
{\it non-symmetric} if $\{y(\R)\}\cap\{-y(\R)\}=\emptyset$.
It was proved in \cite{LLZ} that a prime
characteristic $(\tau, y)$ on $\Sg\in \mathcal{SH}(2n)$ is
symmetric if and only if
$ y(t)=-y(t+\frac{\tau}{2})$ for all $t\in\R$.

There is a long standing conjecture on the number of closed characteristics
on compact convex hypersurfaces in $\R^{2n}$:
\be \,^{\#}\T(\Sg)\ge n, \qquad \forall \; \Sg\in\H(2n). \lb{1.2}\ee

Since the pioneering works \cite{Rab1} of P. Rabinowitz and \cite{Wei1}
of A. Weinstein in 1978 on the existence of at least one closed
characteristic on every hypersurface in $\H(2n)$, the existence of multiple
closed characteristics on $\Sg\in\H(2n)$ has been deeply studied
by many mathematicians. When $n\ge 2$, besides many results under pinching
conditions, in 1987-1988 I. Ekeland-L. Lassoued, I. Ekeland-H. Hofer, and
A, Szulkin (cf. \cite{EkL1}, \cite{EkH1}, \cite{Szu1}) proved
$$ \,^{\#}\T(\Sg)\ge 2, \qquad \forall\,\Sg\in\H(2n). $$
In \cite{HWZ} of 1998, H. Hofer-K. Wysocki-E. Zehnder proved
that $\,^{\#}\T(\Sg)=2$ or $\infty$ holds for every $\Sg\in\H(4)$.
In \cite{LoZ1} of 2002, Y. Long and C. Zhu  proved
$$ \;^{\#}\T(\Sg)\ge [\frac{n}{2}]+1, \qquad \forall\, \Sg\in \H(2n), $$
where we denote by $[a]\equiv\max\{k\in\Z\,|\,k\le a\}$.
In \cite{WHL}, the authors proved the conjecture for $n=3$.
In \cite{LLZ}, the the authors proved the conjecture
for $\Sg\in\mathcal{SH}(2n)$.

Note that in \cite{W2}, the author proved
if $^\#\T(\Sigma)=n$ for some
$\Sigma\in\mathcal{SH}(2n)$ and $n=2$ or $3$, then
any $(\tau, y)\in \T(\Sg)$ is symmetric.
Thus it is natural to conjecture that
\be \,^{\#}\T_s(\Sg)\ge n, \qquad \forall \; \Sg\in\mathcal{SH}(2n), \lb{1.3}\ee
where $T_s(\Sg)$ denotes the set of geometrically distinct symmetric closed
characteristics $(\tau,\, y)$ on $\Sg$.

The following is the main result in this article:

{\bf Theorem 1.1.} {\it We have $ \,^{\#}\T_s(\Sg)\ge 2$ for any
$\Sg\in\mathcal{SH}(2n)$, where $n\ge 2$. }

In this article, let $\N$, $\N_0$, $\Z$, $\Q$, $\R$, and $\C$ denote
the sets of natural integers, non-negative integers, integers, rational
numbers, real numbers, and complex numbers respectively.
Denote by $a\cdot b$ and $|a|$ the standard inner product and norm in
$\R^{2n}$. Denote by $\langle\cdot,\cdot\rangle$ and $\|\cdot\|$
the standard $L^2$-inner product and $L^2$-norm. For an $S^1$-space $X$, we denote
by $X_{S^1}$ the homotopy quotient of $X$ module the $S^1$-action, i.e.,
$X_{S^1}=S^\infty\times_{S^1}X$. We define the functions
\be \left\{\matrix{[a]=\max\{k\in\Z\,|\,k\le a\}, &
E(a)=\min\{k\in\Z\,|\,k\ge a\} , \cr
                   \varphi(a)=E(a)-[a],   \cr}\right. \lb{1.4}\ee
Specially, $\varphi(a)=0$ if $ a\in\Z\,$, and $\varphi(a)=1$ if $a\notin\Z\,$.
In this article we use only $\Q$-coefficients for all homological modules.

\setcounter{equation}{0}
\section{ A variational structure for closed characteristics }

In this section, we transform the problem (\ref{1.1}) into a fixed period
problem of a Hamiltonian system and then study its variational structure.

In the rest of this paper, we fix a $\Sg\in\mathcal{SH}(2n)$ and assume the following
condition on $\Sg$:

\noindent (F) {\bf There exist only finitely many geometrically distinct
symmetric closed characteristics $\{(\tau_j, y_j)\}_{1\le j\le k}$ on $\Sigma$. }

Note that $(\tau, y)\in \mathcal{T}_s(\Sg)$ is a solution of (\ref{1.1}) if
and only if it satisfies the equation
\be
\left\{\matrix{\dot{y}=JN_{\Sigma}(y), \cr
               y(\frac{\tau}{2})=-y(0), \cr }\right. \lb{2.1}\ee
Now we construct a variational structure of closed characteristics
as the following.

{\bf lemma 2.1.} (cf. Proposition 2.2 of \cite{WHL}) {\it  For any sufficiently small $\vartheta\in(0,1)$,
there exists a function $\varphi\equiv \varphi_{\vartheta}\in C^\infty(\R, \;\R^+)$
depending on $\vartheta$ which has $0$ as its unique critical point in $[0, +\infty)$
such that the following hold

(i) $\varphi(0)=0=\varphi^\prime(0)$ and
$ \varphi^{\prime\prime}(0)=1=\lim_{t\rightarrow 0^+}\frac{\varphi^\prime(t)}{t}$.

(ii) $\varphi(t)$ is a polynomial of degree $2$ in a neighborhood of $+\infty$.

(iii) $\frac{d}{dt}\left(\frac{\varphi^\prime(t)}{t}\right)<0$ for $t>0$,
and $\lim_{t\rightarrow +\infty}\frac{\varphi^\prime(t)}{t}<\vartheta$,
i.e., $\frac{\varphi^\prime(t)}{t}$ is strictly decreasing for $t> 0$.

(iv) $\min(\frac{\varphi^\prime(t)}{t}, \varphi^{\prime\prime}(t))\ge \sigma$
for all $t\in \R^+$ and some $\sigma>0$. Consequently, $\varphi$ is strictly
convex on $[0,\,+\infty)$.

(v) In particular, we can choose $\aa\in (1,2)$ sufficiently close to $2$
and $c\in (0,1)$ such that $\vf(t)=ct^{\aa}$ whenever
$\frac{\vf'(t)}{t}\in [\vartheta, 1-\vartheta]$ and $t>0$. }

Let $j: \R^{2n}\rightarrow\R$ be the gauge function of $\Sigma$, i.e.,
$j(\lambda x)=\lambda$ for $x\in\Sigma$ and $\lambda\ge0$, then
$j\in C^3(\R^{2n}\setminus\{0\}, \R)\cap C^0(\R^{2n}, \R)$
and $\Sigma=j^{-1}(1)$. Denote by $\hat{\tau}=\inf_{1\le j\le k}\tau_j$
and $\hat{\sg}=\min\{|y|^2\,|\, y\in\Sigma\}$.

By the same proof of Proposition 2.4 of \cite{WHL}, we have the
following

{\bf Proposition 2.2.} {\it Let $a>\hat{\tau}$,
$\vartheta_a\in\left(0,  \frac{1}{a}\min\{\hat{\tau}, \hat{\sg}\}\right)$
and $\varphi_a$ be a $C^\infty$ function associated to $\vartheta_a$ satisfying
(i)-(iv) of Lemma 2.1. Define the Hamiltonian function $H_a(x)=a\varphi_a(j(x))$
and consider the fixed period problem
\bea
\left\{\matrix{\dot{x}(t)=JH_a^\prime(x(t)) \cr
     x(\frac{1}{2})=-x(0)         \cr }\right. \lb{2.2}\eea
Then the following hold:

(i) $H_a\in C^3(\R^{2n}\setminus\{0\}, \R)\cap C^1(\R^{2n}, \R)$
and there exist $R, r>0$ such that
$$r|\xi|^2\le H^{\prime\prime}_a(x)\xi\cdot\xi\le R|\xi|^2,\quad
\forall x\in \R^{2n}\setminus\{0\},\;\xi\in \R^{2n}.$$

(ii) There exist $\epsilon_1, \epsilon_2\in \left(0, \frac{1}{2}\right)$
and $C\in \R$, such that
$$\frac{\epsilon_1 |x|^2}{2}-C\le H_a(x)\le\frac{\epsilon_2 |x|^2}{2}+C,
\quad \forall x\in\R^{2n}.$$

(iii) Solutions of (\ref{2.2}) are $x\equiv0$ and $x=\rho y(\tau t)$ with
$\frac{\varphi_a^\prime(\rho)}{\rho}=\frac{\tau}{a}$,
where $(\tau, y)$ is a solution of (\ref{2.1}).
In particular, nonzero solutions of (\ref{2.2})
are in one to one correspondence with solutions of (\ref{2.1})
with period $\tau<a$.

(iv) There exists $r_0>0$ independent of $a$ and there
exists $\mu_a>0$ depending on $a$ such that}
$$ H_a^{\prime\prime}(x)\xi\cdot\xi\ge {2ar_0}|\xi|^2,
    \qquad{\rm for}\quad 0<|x|\le \mu_a,\; \xi\in\R^{2n}. $$

In the following, we will use the Clarke-Ekeland dual action principle.
As usual, the Fenchel transform of a function
$F: \R^{2n}\rightarrow\R$ is defined by
\be F^\ast (y)=\sup\{x\cdot y-F(x)\;|\; x\in \R^{2n}\}. \lb{2.3}\ee
Following Proposition 2.2.10 of \cite{Eke3}, Lemma 3.1 of \cite{Eke1}
and the fact that $F_1\le F_2\Leftrightarrow F_1^\ast\ge F_2^\ast$, we
have:

{\bf Proposition 2.3.} {\it Let $H_a$ be a function
defined in Proposition 2.2 and $G_a=H_a^\ast$ the Fenchel
transform of $H_a$. Then we have

(i) $G_a\in C^2(\R^{2n}\setminus\{0\}, \R)\cap C^1(\R^{2n}, \R)$
and $$G_a^\prime(y)=x\Leftrightarrow y=H_a^\prime(x)\Rightarrow
H_a^{\prime\prime}(x)G_a^{\prime\prime}(y)=1.$$

(ii) $G_a$ is strictly convex. Let $R$ and $r$ be the real
numbers given by (i) of Proposition 2.2. Then we have
$$R^{-1}|\xi|^2\le G^{\prime\prime}_a(y)\xi\cdot\xi\le  r^{-1}|\xi|^2,
\quad \forall y\in \R^{2n}\setminus\{0\},\;\xi\in \R^{2n}.$$

(iii) Let $\epsilon_1, \epsilon_2, C$ be the real numbers given by (ii)
of Proposition 2.2. Then we have
$$\frac{|x|^2}{2\epsilon_2}-C\le G_a(x)\le\frac{|x|^2}{2\epsilon_1}+C,
\quad \forall x\in\R^{2n}.$$

(iv) Let $r_0>0$ be the constant given by (iv) of Proposition 2.2.
Then there exists $\eta_a>0$ depending on $a$ such that the following holds
$$G_a^{\prime\prime}(y)\xi\cdot\xi\le \frac{1}{2ar_0}|\xi|^2,
\qquad{\rm for}\quad 0<|y|\le \eta_a,\; \xi\in\R^{2n}. $$

(v) In particular, let $H_a=a\varphi_a(j(x))$ with $\varphi_a$ satisfying
further (v) of Lemma 2.1. Then we have
$G_a(\mu j^\prime(z))=c_1\mu^\beta$ when $z\in \Sigma$
and $\mu j^\prime(z)\in \{H_a^\prime (x)\;|\; H_a(x)=acj(x)^\alpha\}$,
where $c$ is given by (v) of Lemma 2.1, $c_1>0$ is some constant
depending on $a$ and $\alpha^{-1}+\beta^{-1}=1$.} \hfill\hb

Now we apply the dual action principle to problem (\ref{2.3}). Let
\be L^2\left(\R\left/\frac{}{}\right.\left(\frac{1}{2}\Z\right), \R^{2n}\right)=
\{u\in L^2(\R,\R^{2n})|u(t+1/2)=-u(t)\}.
\lb{2.4}\ee
Define a linear operator
$M: L^2\left(\R\left/\frac{}{}\right.\left(\frac{1}{2}\Z\right), \R^{2n}\right)\rightarrow L^2\left(\R\left/\frac{}{}\right.\left(\frac{1}{2}\Z\right), \R^{2n}\right)$ by
\be \frac{d}{dt}Mu(t)=u(t).\lb{2.5}\ee

{\bf Lemma 2.4.} {\it  $M$ is a compact operator from
$L^2\left(\R\left/\frac{}{}\right.\left(\frac{1}{2}\Z\right), \R^{2n}\right)$
into itself and $M^\ast=-M$.}

{\bf Proof.} Note that $M$ sends $L^2\left(\R\left/\frac{}{}\right.\left(\frac{1}{2}\Z\right), \R^{2n}\right)$
into $W^{1, 2}\left([0, 1/2], \R^{2n}\right)$,
and the identity map from $W^{1, 2}\left([0, 1/2], \R^{2n}\right)$
to $L^2\left(\R\left/\frac{}{}\right.\left(\frac{1}{2}\Z\right), \R^{2n}\right)$
is compact by the Rellich-Kondrachov theorem. Hence $M$ is compact.

To check it is anti-symmetric, we use integrate by parts:
\bea \int_0^{1/2}(Mu, v)dt=-\int_0^{1/2}(u, Mv)dt+(Mu, Mv)|_0^{1/2}.
\nn\eea
and the last term vanishes since
$Mu(1/2)=-Mu(0)$ and $Mv(1/2)=-Mv(0)$.
Hence $M$ is anti-symmetric.
\hfill\hb

The dual action functional on $L^2\left(\R\left/\frac{}{}\right.\left(\frac{1}{2}\Z\right), \R^{2n}\right)$ is defined by
\be \Psi_a(u)=\int_0^{1/2}\left(\frac{1}{2}Ju\cdot Mu+G_a(-Ju)\right)dt,
\lb{2.6}\ee
where $G_a$ is given by Proposition 2.3.

{\bf Proposition 2.5.} {\it The functional $\Psi_a$ is bounded from
below on $L^2\left(\R\left/\frac{}{}\right.\left(\frac{1}{2}\Z\right), \R^{2n}\right)$.}

{\bf Proof.} For any $u\in L^2\left(\R\left/\frac{}{}\right.\left(\frac{1}{2}\Z\right), \R^{2n}\right)$,
we represent
$u$ by its Fourier series
\bea u(t)=\sum_{k\in 2\Z+1} e^{2k\pi Jt}x_k, \quad x_k\in\R^{2n}.\lb{2.7}\eea
Then we have
\bea Mu(t)=-J\sum_{k\in 2\Z+1} \frac{1}{2\pi k}e^{2k\pi Jt}x_k.\lb{2.8}\eea
Hence
\bea \frac{1}{2}\langle Ju,\; Mu\rangle=
-\frac{1}{2}\sum_{k\in 2\Z+1} \frac{1}{2\pi k}|x_k|^2
\ge -\frac{1}{4\pi}\|u\|^2.\lb{2.9}
\eea
By (\ref{2.6}), we have
\bea \Psi_a(u)&=&\int_0^{1/2}\left(\frac{1}{2}Ju\cdot Mu+G_a(-Ju)\right)dt\nn\\
&\ge&\frac{1}{2}\langle Ju,\; Mu\rangle
+\int_0^{1/2}\left(\frac{|u|^2}{2\epsilon_2}-C\right)dt.\nn\\
&\ge&\left(\frac{1}{2\epsilon_2}-\frac{1}{4\pi}\right)\|u\|^2-C\nn\\
&\ge& C_1\|u\|^2-C \lb{2.10}
\eea
for some constant $C_1>0$, where in the first inequality, we have
used (iii) of Proposition 2.3. Hence the proposition holds. \hfill\hb

{\bf Proposition 2.6.} {\it The functional $\Psi_a$
is $C^{1, 1}$ on $L^2\left(\R\left/\frac{}{}\right.\left(\frac{1}{2}\Z\right), \R^{2n}\right)$
and satisfies the Palais-Smale condition. Suppose $x$ is a solution of (\ref{2.2}), then $u=\dot{x}$ is a critical
point of $\Psi_a$. Conversely, suppose $u$ is a critical point
of $\Psi_a$, then $Mu$ is a solution of (\ref{2.2}). In particular, solutions
of (\ref{2.2}) are in one to one correspondence with critical
points of $\Psi_a$. }

{\bf Proof.} By (ii) of Proposition 2.3 and the same proof of Proposition 3.3
on p.33 of \cite{Eke1}, we have $\Psi_a$
is $C^{1, 1}$ on $L^2\left(\R\left/\frac{}{}\right.\left(\frac{1}{2}\Z\right), \R^{2n}\right)$.
By (\ref{2.10}) and the proof of Lemma 5.2.8 of \cite{Eke3}, we have $\Psi_a$
 satisfies the Palais-Smale condition.

By (\ref{2.6}), we have
\be \langle\Psi^\prime_a(u), v\rangle=\langle Mu, Jv\rangle-\langle G_a^\prime(-Ju), Jv\rangle,
\lb{2.11}\ee
where we use the fact that
\be Mu(t)=\int_0^tu(s)ds-\frac{1}{2}\int_0^{1/2}u(s)ds,
\lb{2.12}\ee
and  $MJu(t)=JMu(t)$. Hence $\Psi^\prime_a(u)=0$
if and only if $Mu=G_a^\prime(-Ju)$, where we used the fact
$G_a^\prime(-Ju(\frac{1}{2}))=G_a^\prime(Ju(0))=-G_a^\prime(-Ju(0))$.
Taking Frenchel dual we have $-Ju=H_a^\prime(Mu)$, i.e.,
$u=JH_a^\prime(Mu)$. Hence $Mu$ is a solution of (\ref{2.2}).
The converse is obvious. \hfill\hb

{\bf Proposition 2.7.} {\it We have $\Psi_a(u_a)<0$ for every critical
point $u_a\not= 0$ of $\Psi_a$. }

{\bf Proof}. By Propositions 2.2 and 2.6, we have $u_a=\dot{x}_a$ and
$x_a=\rho_ay(\tau t)$ with
\be  \frac{\varphi_a'(\rho_a)}{\rho_a}=\frac{\tau}{a}.  \lb{2.13}\ee
Hence we have
\bea \Psi_a(u_a)
&=& \int_0^{1/2}\left(\frac{1}{2}J\dot x_a\cdot x_a+G_a(-J\dot x_a)\right)dt\nn\\
&=& -\frac{1}{4}\langle H_a'(x_a),\; x_a\rangle+\int_0^{1/2}G_a(H_a'(x_a))dt\nn\\
&=& \frac{1}{4}a\varphi_a'(\rho_a)\rho_a-\frac{1}{2}a\varphi_a(\rho_a). \lb{2.14}
\eea
Here the second equality follows from (\ref{2.2}) and the third
equality follows from (i) of Proposition 2.3 and (\ref{2.3}).

Let $f(t)=\frac{1}{2}a\varphi_a^\prime(t)t-a\varphi_a(t)$ for $t\ge 0$. Then
we have $f(0)=0$ and
$f'(t)=\frac{a}{2}(\varphi_a^{\prime\prime}(t)t-\varphi_a^\prime(t))<0$
since $\frac{d}{dt}(\frac{\varphi_a^\prime(t)}{t})<0$ by (iii) of
Lemma 2.1. This together with (\ref{2.13}) yield the proposition. \hfill\hb

We have a natural $S^1$-action on $L^2\left(\R\left/\frac{}{}\right.\left(\frac{1}{2}\Z\right), \R^{2n}\right)$
 defined by
\bea\theta\ast u(t)=u(\theta+t),\quad\forall\theta\in S^1\equiv\R/\Z, t\in\R.\lb{2.15}\eea
Then we have

{\bf Lemma 2.8.} {\it The functional $\Psi_a$ is $S^1$-invariant. }

{\bf Proof.} Note that we have the following

{\bf Claim. } {\it We have $M(\theta\ast u)=\theta\ast(Mu)$.}

In fact, by (\ref{2.12}), we have
\bea M(\theta\ast u)(t)&=&\int_0^t\theta\ast u(s)ds-\frac{1}{2}\int_0^{1/2}\theta\ast u(s)ds
\nn\\&=&\int_0^t u(\theta+s)ds-\frac{1}{2}\int_0^{1/2} u(\theta+s)ds
\nn\\&=&\int_\theta^{t+\theta}u(s)ds-\frac{1}{2}\int_\theta^{1/2+\theta}\cdot u(s)ds
\nn\eea
On the other hand, we have
\bea \theta\ast (Mu)(t)&=&\theta\ast\left(\int_0^t u(s)ds-\frac{1}{2}\int_0^{1/2} u(s)ds\right)
\nn\\&=&\int_0^{t+\theta}u(s)ds-\frac{1}{2}\int_0^{1/2} u(s)ds
\nn\\&=&\int_0^{\theta} u(s)ds+\int_\theta^{t+\theta} u(s)ds
-\frac{1}{2}\int_0^{\theta} u(s)ds -\frac{1}{2}\int_\theta^{1/2} u(s)ds
\nn\\&=&\frac{1}{2}\int_0^{\theta} u(s)ds+\int_\theta^{t+\theta} u(s)ds
-\frac{1}{2}\int_\theta^{1/2} u(s)ds
\nn\\&=&-\frac{1}{2}\int_{1/2}^{\theta+1/2} u(s)ds+\int_\theta^{t+\theta} u(s)ds
-\frac{1}{2}\int_\theta^{1/2} u(s)ds\lb{2.16}
\\&=&\int_\theta^{t+\theta}u(s)ds-\frac{1}{2}\int_\theta^{1/2+\theta}\cdot u(s)ds,\nn
\eea
where in (\ref{2.16}), we use the fact $u(t+1/2)=-u(t)$.
Hence the claim holds.

Now we have
\bea \Psi_a(\theta\ast u)&=&\int_0^{1/2}\left(\frac{1}{2}J(\theta\ast u)\cdot M(\theta\ast u)+G_a(-J(\theta\ast u))\right)dt,
\nn\\&=&\int_0^{1/2}\left(\frac{1}{2}\theta\ast (Ju)\cdot \theta\ast(Mu)+G_a(\theta\ast (-Ju))\right)dt,
\nn\\&=&\int_\theta^{\theta+1/2}\left(\frac{1}{2}Ju\cdot Mu+G_a(-Ju)\right)dt
\nn\\&=&\int_\theta^{1/2}\left(\frac{1}{2}Ju\cdot Mu+G_a(-Ju)\right)dt
+\int_{1/2}^{\theta+1/2}\left(\frac{1}{2}Ju\cdot Mu+G_a(-Ju)\right)dt
\nn\\&=&\int_\theta^{1/2}\left(\frac{1}{2}Ju\cdot Mu+G_a(-Ju)\right)dt
+\int_{0}^{\theta}\left(\frac{1}{2}(-Ju)\cdot (-Mu)+G_a(Ju)\right)dt
\nn\\&=&\int_\theta^{1/2}\left(\frac{1}{2}Ju\cdot Mu+G_a(-Ju)\right)dt
=\Psi_a(u),\nn
\eea
where in the above computation, we use  $u(t+1/2)=-u(t)$
and $G_a(x)=G_a(-x)$, which follows from $\Sg=-\Sg$. Hence the proposition holds.\hfill\hb

For any $\kappa\in\R$, we denote by
\bea \Lambda_a^\kappa=\left\{u\in L^2\left(\R\left/\frac{}{}\right.\left(\frac{1}{2}\Z\right), \R^{2n}\right)
\;\left|\frac{}{}\right.\;\Psi_a(u)\le\kappa\right\}.\lb{2.17}\eea
For a critical point $u$ of $\Psi_a$, we denote by
\bea \Lambda_a(u)=\Lambda_a^{\Psi_a(u)}
=\left\{w\in L^2\left(\R\left/\frac{}{}\right.\left(\frac{1}{2}\Z\right), \R^{2n}\right)
 \;\left|\frac{}{}\right.\; \Psi_a(w)\le\Psi_a(u)\right\}.\lb{2.18}\eea
Clearly, both sets are $S^1$-invariant. Since the
$S^1$-action preserves $\Psi_a$, if $u$ is a critical
point of $\Psi_a$, then the whole orbit $S^1\cdot u$ is formed by
critical points of $\Psi_a$. Denote by $crit(\Psi_a)$ the set of
critical points of $\Psi_a$. Note that by the condition $(F)$,
(iii) of Proposition 2.2 and Proposition 2.6,
the number of critical orbits of $\Psi_a$ is finite.
Hence as usual we can make the following definition.

{\bf Definition 2.9.} {\it Suppose $u$ is a nonzero critical
point of $\Psi_a$, and $\Nn$ is an $S^1$-invariant
open neighborhood of $S^1\cdot u$ such that
$crit(\Psi_a)\cap(\Lambda_a(u)\cap \Nn)=S^1\cdot u$. Then
the $S^1$-critical modules of $S^1\cdot u$ is defined by
\bea C_{S^1,\; q}(\Psi_a, \;S^1\cdot u)
&=&H_{S^1,\; q}(\Lambda_a(u)\cap\Nn,\;
(\Lambda_a(u)\setminus S^1\cdot u)\cap\Nn)\nn\\
&\equiv&H_{q}((\Lambda_a(u)\cap\Nn)_{S^1},\;
((\Lambda_a(u)\setminus S^1\cdot u)\cap\Nn)_{S^1}),\lb{2.19}
\eea
where $H_{S^1,\;\ast}$ is the $S^1$-equivariant homology with
rational coefficients in the sense of A. Borel
(cf. Chapter IV of \cite{Bor1}).}

By the same argument as Proposition 3.2 of \cite{WHL},
we have the following proposition for critical modules.

{\bf Proposition 2.10.} {\it The critical module $C_{S^1,\;
q}(\Psi_a, \;S^1\cdot u)$ is independent of the choice of $H_a$
defined in Proposition 2.2 in the sense that if $x_i$ are solutions
of (\ref{2.2}) with Hamiltonian functions $H_{a_i}(x)\equiv
a_i\varphi_{a_i}(j(x))$ for $i=1$ and $2$ respectively such that
both $x_1$ and $x_2$ correspond to the same closed characteristic
$(\tau, y)$ on $\Sigma$. Then we have
\be C_{S^1,\; q}(\Psi_{a_1}, \;S^1\cdot\dot {x}_1) \cong
  C_{S^1,\; q}(\Psi_{a_2}, \;S^1\cdot \dot {x}_2), \quad \forall q\in \Z.
\lb{2.20}\ee
In other words, the critical modules are invariant for all $a>\tau$ and
$\varphi_a$ satisfying (i)-(iv) of Lemma 2.1. }

In order to compute the critical modules, as in p.35 of \cite{Eke1}
and p.219 of \cite{Eke3} we introduce the following.

{\bf Definition 2.11.} {\it Suppose $u$ is a nonzero critical
point of $\Psi_a$. Then the formal Hessian of $\Psi_a$ at $u$ is defined by
\be Q_a(v,\; v)=\int_0^{1/2} (Jv\cdot Mv+G_a^{\prime\prime}(-Ju)Jv\cdot Jv)dt,
    \lb{2.21}\ee
which defines an orthogonal splitting $L^2\left(\R\left/\frac{}{}\right.\left(\frac{1}{2}\Z\right), \R^{2n}\right)=E_-\oplus E_0\oplus E_+$ of
$L^2\left(\R\left/\frac{}{}\right.\left(\frac{1}{2}\Z\right), \R^{2n}\right)$ into negative, zero and positive subspaces. The index
of $u$ is defined by $i(u)=\dim E_-$ and the nullity of $u$ is defined by
$\nu(u)=\dim E_0$. }

Next we show that the index and nullity defined as above are the Morse
index and nullity of a corresponding functional on a finite dimensional
subspace of $L^2\left(\R\left/\frac{}{}\right.\left(\frac{1}{2}\Z\right), \R^{2n}\right)$.

{\bf Lemma 2.12.} {\it Let  $\Psi_a$ be a functionals defined by (\ref{2.6}).
Then there exists a finite dimensional $S^1$-invariant subspace $X$  of
$L^2\left(\R\left/\frac{}{}\right.\left(\frac{1}{2}\Z\right), \R^{2n}\right)$
 and a  $S^1$-equivariant map
$h_a: X\rightarrow X^\perp$ such that the following hold

(i) For $g\in X$, each function $h\mapsto\Psi_a(g+h)$
has $h_a(g)$ as the unique minimum in $X^\perp$.

Let $\psi_a(g)=\Psi_a(g+h_a(g))$. Then we have

(ii) The function $\psi_a$ is $C^1$ on $X$ and $S^1$-invariant.
$g_a$ is a critical point of $\psi_a$ if and only if $g_a+h_a(g_a)$
is a critical point of $\Psi_a$.

(iii) If $g_a\in X$ and $H_a$ is $C^k$ with $k\ge 2$ in a neighborhood
of the trajectory of $g_a+h_a(g_a)$, then  $\psi_a$ is $C^{k-1}$ in a
neighborhood of $g_a$. In particular,
if $g_a$ is a nonzero critical point of $\psi_a$,
then  $\psi_a$ is $C^2$ in a neighborhood of the critical
orbit $S^1\cdot g_a$.
The index and nullity of $\Psi_a$ at $g_a+h_a(g_a)$ defined
in Definition 2.11 coincide with the Morse index and nullity of
$\psi_a$ at $g_a$.

(iv) For any $\kappa\in\R$, we denote by
\bea \widetilde{\Lambda}_a^\kappa=\{g\in X \;|\; \psi_a(g)\le\kappa\}. \lb{2.22}\eea
Then the natural embedding
$\widetilde{\Lambda}_a^\kappa \hookrightarrow  {\Lambda}_a^\kappa $
given by $g\mapsto g+h_a(g)$
is an $S^1$-equivariant homotopy equivalence.
}

{\bf Proof.} By (ii) of Proposition 2.3, we have
\bea (G_a^\prime(u)-G_a^\prime(v), u-v)\ge
\omega|u-v|^2,\quad \forall u, v\in\R^{2n},
\lb{2.23}\eea
for some $\omega>0$. Hence we can use the proof of
Proposition 3.9 of \cite{Vit1} to obtain $X$ and $h_a$.
In fact, $X$ is the subspace of $L^2\left(\R\left/\frac{}{}\right.\left(\frac{1}{2}\Z\right), \R^{2n}\right)$
generated by the eigenvectors of $-JM$ whose eigenvalues are less than
$-\frac{\omega}{2}$ and $h_a(g)$ is defined by the equation
\bea \frac{\partial}{\partial h}\Psi_a(g+h_a(g))=0,\nn\eea
then (i)-(iii) follows from Proposition 3.9 of \cite{Vit1}.
(iv) follows from Lemma 5.1 of \cite{Vit1}.\hfill\hb

Note that $\Psi_a$ is not $C^2$ in general, and then we can not apply
Morse theory to $\Psi_a$ directly. After the finite dimensional
approximation, the function $\psi_a$ has much better differentiability,
which allows us to apply the Morse theory to study its property.

{\bf Proposition 2.13.} {\it Let $\Psi_a$ be a
functional defined by (\ref{2.6}), and $u_a=\dot{x}_a$ be the critical
point of $\Psi_a$ so that $x_a$ corresponds to a closed
characteristic $(\tau,y)$ on $\Sigma$. Then the
nullity $\nu(u_a)$ of the functional $\Psi_a$ at its
critical point $u_a$ is the number of linearly independent
solutions of the boundary value problem
\bea
\left\{\matrix{\dot{\xi}(t)=JH_a^{\prime\prime}(x_a(t))\xi \cr
     \xi(\frac{1}{2})=-\xi(0)         \cr }\right. \lb{2.24}\eea}

{\bf Proof.} By (\ref{2.21}), we have
\bea Q_a(v,\; w)&=&\int_0^{1/2} (Jv\cdot Mw+G_a^{\prime\prime}(-Ju)Jv\cdot Jw)dt,
\nn\\&=&\langle Mw, Jv\rangle+\langle(H_a^{\prime\prime}(x_a(t))^{-1}Jw, Jv\rangle
    \lb{2.25}\eea
where we have used (\ref{2.2}) and (i) of Proposition 2.3.
Now $w\in E_0$ if and only if $Q_a(v,\; w)=0$ for any
$v\in L^2\left(\R\left/\frac{}{}\right.\left(\frac{1}{2}\Z\right), \R^{2n}\right)$.
Hence we must have $Mw+(H_a^{\prime\prime}(x_a(t))^{-1}Jw=0$,
i.e., we have $w=JH_a^{\prime\prime}(x_a(t))Mw$.
Hence $Mw$ solves (\ref{2.25}).\hfill\hb

Denote by $R(t)$ the fundamental solution
of the linearized system
\be \dot \xi(t)=JH_a^{\prime\prime}(x_a(t))\xi(t), \lb{2.26}\ee
Then we have the following

{\bf Proposition 2.14.} {\it In an appropriate coordinates there holds
\bea R(1/2)=\left(\matrix{A\quad B  \cr
                        0\quad C   \cr }\right)
\quad {\rm with}\quad A=\left(\matrix{-1\quad -\gamma  \cr
                        0\quad -1   \cr }\right),\nn
\eea
with $\gamma>0$ and $C$ is independent of $H_a$. }

{\bf Proof.} Note that by Lemma 1.6.11 of \cite{Eke3}, we have
\be R(t)T_{y(0)}\Sigma\subset T_{y(\tau t)}\Sigma.\lb{2.27}\ee
Differentiating (\ref{2.2}) and use the fact $x_a(t+1/2)=-x_a(t)$,
we have
\be R(1/2)\dot x_a(0)=-\dot x_a(0).\lb{2.28}\ee
Let
\be x_a(\rho, t)=\rho y\left(\frac{\tau t}{T_\rho}\right) \quad {\rm with}\;
\frac{\tau}{T_\rho}=\frac{a\varphi_a^\prime(\rho)}{\rho}.\lb{2.29}\ee
Then we have $x_a(\rho, T_\rho/2)=-x_a(\rho, 0)$. Differentiating
it with respect to $\rho$ and using (\ref{2.29}) together with
$\dot x_a(1/2)=-\dot x_a(0)$, we get
$$
-\frac{\tau}{2a}\frac{d}{d\rho}\left(\frac{\rho}{\varphi_a^\prime(\rho)}\right)\dot x_a(0)
+R(1/2)\rho^{-1}x_a(0)=-\rho^{-1}x_a(0).$$
Hence we have
\be R(1/2)x_a(0)=-x_a(0)+
\frac{\rho\tau}{2a}\frac{d}{d\rho}\left(\frac{\rho}{\varphi_a^\prime(\rho)}\right)\dot x_a(0)
=x_a(0)+\gamma \dot x_a(0),\lb{2.30}\ee
where $\gamma>0$ since $\frac{d}{d\rho}\left(\frac{\rho}{\varphi_a^\prime(\rho)}\right)>0$
by (iii) of Proposition 2.1. For any $w\in\R^{2n}$, we have
\bea H_a^{\prime\prime}(x_a)w
&=& a\varphi_a^{\prime\prime}(j(x_a))(j^\prime(x_a), w)j^\prime(x_a)
+a\varphi_a^{\prime}(j(x_a))j^{\prime\prime}(x_a)w\nn\\
&=& a\varphi_a^{\prime\prime}(j(x_a))(j^\prime(y), w)j^\prime(y)
+\tau j^{\prime\prime}(y)w. \lb{2.31}
\eea
The last equality follows from (iii) of Proposition 2.2.
Let $z(t)=R(t)z(0)$ for $z(0)\in T_{y(0)}\Sigma$.
Then by (\ref{2.27}), we have $\dot z(t)=\tau j^{\prime\prime}(y(t))z(t)$.
Therefore $R(1/2)|_{T_{y(0)}}\Sigma$ is
independent of the choice of $H_a$ in Proposition 2.2. Summing up,
we have proved that in an appropriate coordinates there holds
\bea R(1/2)=\left(\matrix{A\quad B  \cr
                        0\quad C   \cr }\right)
\quad {\rm with}\quad A=\left(\matrix{-1\quad -\gamma  \cr
                        0\quad -1   \cr }\right),\nn
\eea
with  $C$ is independent of $H_a$,
where we use $\{-\dot x_a(0), x_a(0), e_1,\ldots,e_{2n-2}\}$
as an basis of $\R^{2n}$. \hfill\hb

{\bf Proposition 2.15.} {\it Let $\Psi_a$ be a
functional defined by (\ref{2.6}), and $u$ be a nonzero critical
point of $\Psi_a$. Then we have
\be  C_{S^1,\; q}(\Psi_a, \;S^1\cdot u)=0,\qquad \forall q\notin[i(u), i(u)+\nu(u)-1].
\lb{2.32}\ee}

{\bf Proof.} By (iv) of Lemma 2.12, we have
\be  C_{S^1,\; q}(\Psi_a, \;S^1\cdot u)
\simeq C_{S^1,\; q}(\psi_a, \;S^1\cdot u),\lb{2.33}\ee
where $C_{S^1,\; q}(\psi_a, \;S^1\cdot u)=H_{S^1,\; q}(\widetilde{\Lambda}_a(u)\cap\Nn,\;
(\widetilde{\Lambda}_a(u)\setminus S^1\cdot u)\cap\Nn)$
and $\Nn$ is an $S^1$-invariant
open neighborhood of $S^1\cdot u$ such that
$crit(\psi_a)\cap(\widetilde{\Lambda}_a(u)\cap \Nn)=S^1\cdot u$.
By (iii) of Lemma 2.12, the functional $\psi_a$ is $C^2$ near $S^1\cdot u$.
Thus we can use  the Gromoll-Meyer theory in the equivariant sense
to obtain the proposition.\hfill\hb

Recall that for a principal $U(1)$-bundle $E\to B$, the Fadell-Rabinowitz index
(cf. \cite{FaR1}) of $E$ is defined to be $\sup\{k\;|\, c_1(E)^{k-1}\not= 0\}$,
where $c_1(E)\in H^2(B,\Q)$ is the first rational Chern class. For a $U(1)$-space,
i.e., a topological space $X$ with a $U(1)$-action, the Fadell-Rabinowitz index is
defined to be the index of the bundle $X\times S^{\infty}\to X\times_{U(1)}S^{\infty}$,
where $S^{\infty}\to CP^{\infty}$ is the universal $U(1)$-bundle.
For any $\kappa\in\R$, we denote by
\be \Psi_a^{\kappa-}=\left\{w\in L^2\left(\R\left/\frac{}{}\right.\left(\frac{1}{2}\Z\right), \R^{2n}\right)
 \;\left|\frac{}{}\right.\; \Psi_a(w)<\kappa\right\}. \lb{2.34}\ee
Then as in P.218 of \cite{Eke3}, we define
\be c_i=\inf\{\delta\in\R\;|\: \hat I(\Psi_a^{\kappa-})\ge i\},\lb{2.35}\ee
where $\hat I$ is the Fadell-Rabinowitz index given above. Then as Proposition 3
in P.218 of \cite{Eke3}, we have

{\bf Proposition 2.16.} {\it Every $c_i$ is a critical value of $\Psi_a$. If
$c_i=c_j$ for some $i<j$, then there are infinitely many geometrically
distinct symmetric closed characteristics on $\Sg$.}

By a similar argument as Proposition 3.5  of \cite{W1}
and Proposition 2.15, we have

{\bf Proposition 2.17.} {\it Suppose $u$ is the critical point of $\Psi_a$ found
in Proposition 2.16. Then we have
\be \Psi_a(u)=c_i,\qquad C_{S^1,\; 2(i-1)}(\Psi_a, \;S^1\cdot u)\neq 0. \lb{2.36}\ee
In particular, we have $i(u)\le 2(i-1)\le i(u)+\nu(u)-1$.}

\setcounter{equation}{0}
\section{Index iteration theory for symmetric closed characteristics}

In this section, we study the index iteration theory for symmetric closed characteristics.

Note that if $(\tau, y)\in \mathcal{T}_s(\Sg)$,
then $((2m-1)\tau, y)$ is a solution of (\ref{2.1})  for any $m\in\N$.
Thus $((2m-1)\tau, y)$ corresponds to a critical point
of $\Psi_a$ via Propositions 2.2 and 2.6, we denote
it by $u^{2m-1}$. First note that we have the following

{\bf Lemma 3.1.} {\it Suppose $u^{2m-1}$ is a nonzero critical
point of $\Psi_a$ such that $u$ corresponds to $(\tau, y)\in \mathcal{T}_s(\Sg)$.
Let $H(x)=j(x)^2$, where $j$ is the gauge function of $\Sg$.
Then $i(u^{2m-1})$ equals the index of the following quadratic form
\be Q_{(2m-1)\tau/2}(\xi,\; \xi)=\int_0^{(2m-1)\tau/2} (J\dot \xi\cdot \xi+(H^{\prime\prime}(y(t)))^{-1}J\dot\xi\cdot J\dot\xi)dt,
    \lb{3.1}\ee
where $\xi\in W^{1, 2}\left(\R\left/\frac{}{}\right.\left(\frac{(2m-1)\tau}{2}\Z\right), \R^{2n}\right)
\equiv\{w\in W^{1, 2}(\R,\R^{2n})|w(t+\frac{(2m-1)\tau}{2})=-w(t)\}.$.
Moreover, we have $\nu(u^{2m-1})={\rm nullity}Q_{(2m-1)\tau/2}-1$.
}

{\bf Proof.} By a similar argument as in proposition 1.7.5
and P.36 of \cite{Eke3} and Proposition 3.5 of \cite{WHL},
we obtain the lemma. \hfill\hb

 Suppose $u^{2k-1}$ is a nonzero critical
point of $\Psi_a$ such that $u$ corresponds to $(\tau, y)\in \mathcal{T}_s(\Sg)$.
Then for any $\omega\in \U$, let
\be Q^\omega_{(2k-1)\tau/2}(\xi,\; \xi)=\int_0^{(2k-1)\tau/2} (J\dot \xi\cdot \xi+(H^{\prime\prime}(y(t)))^{-1}J\dot\xi\cdot J\dot\xi)dt,
    \lb{3.2}\ee
where $\xi\in E^\omega_{(2k-1)\tau/2}\equiv
\{u\in W^{1, 2}([0, (2k-1)\tau/2],\C^{2n})|w(\frac{(2k-1)\tau}{2})=\omega w(0)\}.$.

Clearly the quadratic form $Q_{(2m-1)\tau/2}$ on the real Hilbert space
$W^{1, 2}\left(\R\left/\frac{}{}\right.\left(\frac{(2m-1)\tau}{2}\Z\right), \R^{2n}\right)$
and the Hermitian form $Q^{-1}_{(2m-1)\tau/2}$ on the complex Hilbert space
$E^{-1}_{(2m-1)\tau/2}$ have the same index.

If $\omega^{2m-1}=-1$, we identify $E^\omega_{\tau/2}$
with a subspace of $E^{-1}_{(2m-1)\tau/2}$ via
\be E^\omega_{\tau/2}=
\{u\in W^{1, 2}(\R,\C^{2n})|w(t+\tau/2)=\omega w(t)\}.
\lb{3.3}\ee
Note that if $\xi\in E^\omega_{\tau/2}$, we have
\bea Q^\omega_{(2m-1)\tau/2}(\xi,\; \xi)&=&\int_0^{(2m-1)\tau/2} (J\dot \xi\cdot \xi+(H^{\prime\prime}(y(t)))^{-1}J\dot\xi\cdot J\dot\xi)dt
\nn\\&=&\sum_{k=0}^{2m-1}(\omega\overline{\omega})^k\int_0^{\tau/2}(J\dot \xi\cdot \xi+(H^{\prime\prime}(y(t)))^{-1}J\dot\xi\cdot J\dot\xi)dt
\nn\\&=&(2m-1)Q^\omega_{\tau/2}(\xi, \xi).
    \lb{3.4}\eea

{\bf Lemma 3.2.} {\it  The spaces $E^\omega_{\tau/2}$ for $\omega^{2m-1}=-1$
are orthogonal subspaces of $E^{-1}_{(2m-1)\tau/2}$,
both for the standard Hilbert structure and for $Q^{-1}_{(2m-1)\tau/2}$,
and we have the decomposition
\be E^{-1}_{(2m-1)\tau/2}=\bigoplus_{\omega^{2m-1}=-1}E^\omega_{\tau/2}.
\lb{3.5}\ee
}

{\bf Proof.} Any $\xi\in E^{-1}_{(2m-1)\tau/2}$ can be written as
\bea \xi(t)=\sum_{p\in 2\Z+1}x_p\exp\left(\frac{2i\pi pt}{(2m-1)\tau}\right)
\lb{3.6}\eea
for $q=1, 3,\ldots,4m-3$, denote by $C(q)$ the set of
all $p$ such that $p-q\in (4m-2)\Z$. Thus we may write
\bea \xi(t)=\sum_{q\in 2\Z+1\atop 1\le q\le 4m-3}\xi_q(t), \qquad
\xi_q(t)=\sum_{C(q)}x_p\exp\left(\frac{2i\pi pt}{(2m-1)\tau}\right)
\lb{3.7}\eea
Then we have
\bea \xi_q(t+\tau/2)&=&\sum_{C(q)}x_p\exp\left(\frac{2i\pi pt}{(2m-1)\tau}+\frac{i\pi p}{2m-1}\right)
\nn\\&=&\exp\left(\frac{i\pi q}{2m-1}\right)\xi_q(t).
\lb{3.8}\eea
Thus $\xi_q\in E^\omega_{\tau/2}$ with $\omega=\exp\left(\frac{i\pi q}{2m-1}\right)$,
when $q$ runs from $1, 3,\ldots, 4m-3$, then $\omega$
runs through the $2m-1$ roots of $-1$.

For $\xi\in E^\omega_{\tau/2}$  and  $\eta\in E^\lambda_{\tau/2}$ with $\omega\neq \lambda$
are $2m-1$ roots of $-1$, we have
\bea Q^{-1}_{(2m-1)\tau/2}(\xi,\; \eta)&=&\int_0^{(2m-1)\tau/2} (J\dot \xi\cdot \eta+(H^{\prime\prime}(y(t)))^{-1}J\dot\xi\cdot J\dot\eta)dt
\nn\\&=&\sum_{k=0}^{2m-1}(\omega\overline{\lambda})^k\int_0^{\tau/2}(J\dot \xi\cdot \eta+(H^{\prime\prime}(y(t)))^{-1}J\dot\xi\cdot J\dot\eta)dt
\nn\\&=&0.
    \lb{3.9}\eea
Thus the lemma holds. \hfill\hb

{\bf Definition 3.3.} {\it We define the Bott maps $j_{\tau/2}$ and
$n_{\tau/a}$ from $\U$ to $\Z$ by
\be j_{\tau/2}(\omega)={\rm index} Q^\omega_{\tau/2},
\qquad n_{\tau/2}(\omega)={\rm nullity} Q^\omega_{\tau/2},
\lb{3.10}\ee}
By Lemmas 3.1 and 3.2, we have

{\bf Proposition 3.4.} {\it Suppose $u^{2m-1}$ is a nonzero critical
point of $\Psi_a$ such that $u$ corresponds to $(\tau, y)\in \mathcal{T}_s(\Sg)$.
Then we have
\be i(u^{2m-1})=\sum_{\omega^{2m-1}=-1}j_{\tau/2}(\omega)
\qquad \nu(u^{2m-1})=\sum_{\omega^{2m-1}=-1}n_{\tau/2}(\omega)-1.
\lb{3.11}\ee
}

Note that $j_{\tau/2}(\omega)$ coincide with the function
defined in Definition 1.5.3 of \cite{Eke3} for the linear
Hamiltonian system
\bea
\left\{\matrix{\dot{\xi}(t)=JA(t)\xi \cr
     A(t+\tau/2)=A(t)        \cr }\right. \lb{3.12}\eea
where $A(t)=H^{\prime\prime}(y(t))$. Denote by $i^E(A, k)$ and $\nu^E(A, k)$
the index and nullity of the $k$-th iteration of the system (\ref{3.12})
defined by Ekeland in \cite{Eke3}.
Denote by $i(A, k)$ and $\nu(A, k)$ the Maslov-type index and nullity of
the $k$-th iteration of the system (\ref{3.12}) defined by Conley, Zehnder and Long (cf. \S 5.4 of \cite{Lon4}).
Then we have

{\bf Theorem 3.5.} (cf. Theorem 15.1.1 of \cite{Lon4}) {\it We have
\be i^E(A, k)= i(A, k)-n,\quad \nu^E(A, k)=\nu(A, k),
      \lb{3.13}\ee
for any $k\in\N$.}

{\bf Theorem 3.6.} {\it Suppose $u^{2m-1}$ is a nonzero critical
point of $\Psi_a$ such that $u$ corresponds to $(\tau, y)\in \mathcal{T}_s(\Sg)$.
Then we have
\be i(u^{2m-1})=i_{-1}(A, 2m-1),
\qquad \nu(u^{2m-1})=\nu_{-1}(A, 2m-1)-1.
\lb{3.14}\ee
where $i_{-1}(A, k)$ and $\nu_{-1}(A, k)$
are the Maslov-type index and nullity introduced in
\cite{Lon2}.}

{\bf Proof.} By Corollary 1.5.4 of \cite{Eke3}
and Theorem 9.2.1 of \cite{Lon4} respectively, we have
\bea i^E(A, 4m-2)&=&i^E(A, 2m-1)+i_{-1}^E(A, 2m-1),\nn\\
i(A, 4m-2)&=&i(A, 2m-1)+i_{-1}(A, 2m-1)
\lb{3.15}\eea
and by Lemma 3.1, we have $i(u^{2m-1})=i^E_{-1}(A, 2m-1)$.
Thus the theorem follows from Theorem 3.5.\hfill\hb

Now we compute $i(u^{2m-1})$ via the index iteration method
in \cite{Lon4}. First we recall briefly an index theory for symplectic paths.
All the details can be found in \cite{Lon4}.

As usual, the symplectic group $\Sp(2n)$ is defined by
$$ \Sp(2n) = \{M\in {\rm GL}(2n,\R)\,|\,M^TJM=J\}, $$
whose topology is induced from that of $\R^{4n^2}$. For $\tau>0$ we are interested
in paths in $\Sp(2n)$:
$$ \P_{\tau}(2n) = \{\ga\in C([0,\tau],\Sp(2n))\,|\,\ga(0)=I_{2n}\}, $$
which is equipped with the topology induced from that of $\Sp(2n)$. The
following real function was introduced in \cite{Lon2}:
$$ D_{\om}(M) = (-1)^{n-1}\ol{\om}^n\det(M-\om I_{2n}), \qquad
          \forall \om\in\U,\, M\in\Sp(2n). $$
Thus for any $\om\in\U$ the following codimension $1$ hypersurface in $\Sp(2n)$ is
defined in \cite{Lon2}:
$$ \Sp(2n)_{\om}^0 = \{M\in\Sp(2n)\,|\, D_{\om}(M)=0\}.  $$
For any $M\in \Sp(2n)_{\om}^0$, we define a co-orientation of $\Sp(2n)_{\om}^0$
at $M$ by the positive direction $\frac{d}{dt}Me^{t\ep J}|_{t=0}$ of
the path $Me^{t\ep J}$ with $0\le t\le 1$ and $\ep>0$ being sufficiently
small. Let
\bea
\Sp(2n)_{\om}^{\ast} &=& \Sp(2n)\bs \Sp(2n)_{\om}^0,   \nn\\
\P_{\tau,\om}^{\ast}(2n) &=&
      \{\ga\in\P_{\tau}(2n)\,|\,\ga(\tau)\in\Sp(2n)_{\om}^{\ast}\}, \nn\\
\P_{\tau,\om}^0(2n) &=& \P_{\tau}(2n)\bs  \P_{\tau,\om}^{\ast}(2n).  \nn\eea
For any two continuous arcs $\xi$ and $\eta:[0,\tau]\to\Sp(2n)$ with
$\xi(\tau)=\eta(0)$, it is defined as usual:
$$ \eta\ast\xi(t) = \left\{\matrix{
            \xi(2t), & \quad {\rm if}\;0\le t\le \tau/2, \cr
            \eta(2t-\tau), & \quad {\rm if}\; \tau/2\le t\le \tau. \cr}\right. $$
Given any two $2m_k\times 2m_k$ matrices of square block form
$M_k=\left(\matrix{A_k&B_k\cr
                                C_k&D_k\cr}\right)$ with $k=1, 2$,
as in \cite{Lon4}, the $\;\dm$-product of $M_1$ and $M_2$ is defined by
the following $2(m_1+m_2)\times 2(m_1+m_2)$ matrix $M_1\dm M_2$:
$$ M_1\dm M_2=\left(\matrix{A_1&  0&B_1&  0\cr
                               0&A_2&  0&B_2\cr
                             C_1&  0&D_1&  0\cr
                               0&C_2&  0&D_2\cr}\right). \nn$$  
Denote by $M^{\dm k}$ the $k$-fold $\dm$-product $M\dm\cdots\dm M$. Note
that the $\dm$-product of any two symplectic matrices is symplectic. For any two
paths $\ga_j\in\P_{\tau}(2n_j)$ with $j=0$ and $1$, let
$\ga_0\dm\ga_1(t)= \ga_0(t)\dm\ga_1(t)$ for all $t\in [0,\tau]$.

A special path $\xi_n\in\P_{\tau}(2n)$ is defined by
\be \xi_n(t) = \left(\matrix{2-\frac{t}{\tau} & 0 \cr
                                             0 &  (2-\frac{t}{\tau})^{-1}\cr}\right)^{\dm n}
         \qquad {\rm for}\;0\le t\le \tau.  \lb{3.16}\ee
{\bf Definition 3.7.} (cf. \cite{Lon2}, \cite{Lon4}) {\it For any $\om\in\U$ and
$M\in \Sp(2n)$, define
\be  \nu_{\om}(M)=\dim_{\C}\ker_{\C}(M - \om I_{2n}).  \lb{3.17}\ee
For any $\tau>0$ and $\ga\in \P_{\tau}(2n)$, define
\be  \nu_{\om}(\ga)= \nu_{\om}(\ga(\tau)).  \lb{3.18}\ee

If $\ga\in\P_{\tau,\om}^{\ast}(2n)$, define
\be i_{\om}(\ga) = [\Sp(2n)_{\om}^0: \ga\ast\xi_n],  \lb{3.19}\ee
where the right hand side of (\ref{3.19}) is the usual homotopy intersection
number, and the orientation of $\ga\ast\xi_n$ is its positive time direction under
homotopy with fixed end points.

If $\ga\in\P_{\tau,\om}^0(2n)$, we let $\mathcal{F}(\ga)$
be the set of all open neighborhoods of $\ga$ in $\P_{\tau}(2n)$, and define
\be i_{\om}(\ga) = \sup_{U\in\mathcal{F}(\ga)}\inf\{i_{\om}(\beta)\,|\,
                       \beta\in U\cap\P_{\tau,\om}^{\ast}(2n)\}.
               \lb{3.20}\ee
Then
$$ (i_{\om}(\ga), \nu_{\om}(\ga)) \in \Z\times \{0,1,\ldots,2n\}, $$
is called the index function of $\ga$ at $\om$. }

For any $M\in\Sp(2n)$ and $\om\in\U$, the {\it splitting numbers} $S_M^{\pm}(\om)$
of $M$ at $\om$ are defined by
\be S_M^{\pm}(\om)
     = \lim_{\ep\to 0^+}i_{\om\exp(\pm\sqrt{-1}\ep)}(\ga) - i_{\om}(\ga),
   \lb{3.21}\ee
for any path $\ga\in\P_{\tau}(2n)$ satisfying $\ga(\tau)=M$.

Let $\Omega^0(M)$ be the path connected component
 containing $M=\gamma(\tau)$ of the set
\begin{eqnarray}
  \Omega(M)=\{N\in{\rm Sp}(2n)\,&|&\,\sigma(N)\cap{\bf U}=\sigma(M)\cap{\bf U}\;
{\rm and}\;  \nonumber\\
 &&\qquad \nu_{\lambda}(N)=\nu_{\lambda}(M)\;\forall\,
\lambda\in\sigma(M)\cap{\bf U}\}. \label{3.22}
\end{eqnarray}
Here $\Omega^0(M)$ is called the {\it homotopy component} of $M$ in
${\rm Sp}(2n)$.

In \cite{Lon2}-\cite{Lon4}, the following symplectic matrices were introduced
as {\it basic normal forms}:
\begin{eqnarray}
D(\lambda)=\left(\matrix{\lm & 0\cr
         0  & \lm^{-1}\cr}\right), &\quad& \lm=\pm 2,\lb{3.23}\\
N_1(\lm,b) = \left(\matrix{\lm & b\cr
         0  & \lm\cr}\right), &\quad& \lm=\pm 1, b=\pm1, 0, \lb{3.24}\\
R(\th)=\left(\matrix{\cos\th & -\sin\th\cr
        \sin\th  & \cos\th\cr}\right), &\quad& \th\in (0,\pi)\cup(\pi,2\pi),
                     \lb{3.25}\\
N_2(\om,b)= \left(\matrix{R(\th) & b\cr
              0 & R(\th)\cr}\right), &\quad& \th\in (0,\pi)\cup(\pi,2\pi),
                     \lb{3.26}\end{eqnarray}
where $b=\left(\matrix{b_1 & b_2\cr
               b_3 & b_4\cr}\right)$ with  $b_i\in\R$ and  $b_2\not=b_3$.

Splitting numbers possess the following properties:

{\bf Lemma 3.8.} (cf. \cite{Lon2} and Lemma 9.1.5 of \cite{Lon4}) {\it Splitting
numbers $S_M^{\pm}(\om)$ are well defined, i.e., they are independent of the choice
of the path $\ga\in\P_\tau(2n)$ satisfying $\ga(\tau)=M$ appeared in (\ref{3.21}).
For $\om\in\U$ and $M\in\Sp(2n)$, splitting numbers $S_N^{\pm}(\om)$ are constant
for all $N\in\Om^0(M)$. Moreover, we have
\bea &&S_M^{\pm}(\omega) = 0, \qquad {\it if}\;\;\om\not\in\sg(M).
 \nn\\&&S_M^+(\omega)=S_M^-(\overline{\omega}), \qquad \forall \omega\in\U.
\nn \eea }

{\bf Lemma 3.9.} (cf. \cite{Lon2}, Lemma 9.1.5 of \cite{Lon4})
{\it For any $M_i\in\Sp(2n_i)$ with $i=0$ and $1$, there holds }
\be S^{\pm}_{M_0\dm M_1}(\om) = S^{\pm}_{M_0}(\om) + S^{\pm}_{M_1}(\om),
    \qquad \forall\;\om\in\U. \lb{3.27}\ee

We have the following

{\bf Theorem 3.10.} (cf. \cite{Lon3} and Theorem 1.8.10 of \cite{Lon4}) {\it For
any $M\in\Sp(2n)$, there is a path $f:[0,1]\to\Om^0(M)$ such that $f(0)=M$ and
\be f(1) = M_1\dm\cdots\dm M_l,  \lb{3.28}\ee
where each $M_i$ is a basic normal form listed in (\ref{3.23})-(\ref{3.26})
for $1\leq i\leq l$.}

Now we deduce the index iteration formula for each case in (\ref{3.23})-(\ref{3.26}),
Note that the splitting numbers are computed in List 9.1.12 of \cite{Lon4}.

{\bf Case 1.} {\it $M$ is conjugate to a matrix
$\left(\matrix{ 1 & b\cr
                  0 & 1\cr}\right)$ for some $b>0$.}

In this case, we have $(S_M^+(1), S_M^-(1))=(1, 1)$.
Thus by Theorem 9.2.1 of \cite{Lon4}, we have
\bea i_{-1}(\gamma^{2m-1})&=&\sum_{\omega^{2m-1}=-1}i_\omega(\gamma)
=\sum_{k=1}^{2m-1}i_{\frac{(2k-1)\pi}{2m-1}}(\gamma)
=(2m-1)(i_1(\gamma)+1),\nn\\
\nu_{-1}(\gamma^{2m-1})&=&0.\lb{3.29}\eea

{\bf Case 2.} {\it $M=I_2$, the $2\times 2$ identity matrix}.

In this case, we have $(S_M^+(1), S_M^-(1))=(1, 1)$.
Thus as in Case 1, we have
\bea i_{-1}(\gamma^{2m-1})=(2m-1)(i_1(\gamma)+1),\quad
\nu_{-1}(\gamma^{2m-1})=0.\lb{3.30}\eea

{\bf Case 3.} {\it $M$ is conjugate to a matrix
$\left(\matrix{ 1 & b\cr
                  0 & 1\cr}\right)$ for some $b<0$.}

In this case, we have $(S_M^+(1), S_M^-(1))=(0, 0)$.
Thus by Theorem 9.2.1 of \cite{Lon4}, we have
\bea i_{-1}(\gamma^{2m-1})&=&\sum_{\omega^{2m-1}=-1}i_\omega(\gamma)
=\sum_{k=1}^{2m-1}i_{\frac{(2k-1)\pi}{2m-1}}(\gamma)
=(2m-1)i_1(\gamma),\nn\\
\nu_{-1}(\gamma^{2m-1})&=&0.\lb{3.31}\eea

{\bf Case 4.} {\it $M$ is conjugate to a matrix
$\left(\matrix{ -1 & b\cr
                  0 & -1\cr}\right)$ for some $b<0$.}

In this case, we have $(S_M^+(-1), S_M^-(-1))=(1, 1)$.
Thus by Theorem 9.2.1 of \cite{Lon4}, we have
\bea i_{-1}(\gamma^{2m-1})&=&\sum_{\omega^{2m-1}=-1}i_\omega(\gamma)
=\sum_{k=1}^{2m-1}i_{\frac{(2k-1)\pi}{2m-1}}(\gamma)
\nn\\&=&\sum_{k=1}^{m-1}i_{\frac{(2k-1)\pi}{2m-1}}(\gamma)
+i_{-1}(\gamma)+\sum_{k=m+1}^{2m-1}i_{\frac{(2k-1)\pi}{2m-1}}(\gamma)
\nn\\&=&(m-1)i_1(\gamma)+i_1(\gamma)-1+(m-1)(i_1(\gamma)-1+1)
\nn\\&=&(2m-1)i_1(\gamma)-1,
\nn\\
\nu_{-1}(\gamma^{2m-1})&=&1.\lb{3.32}\eea

{\bf Case 5.} {\it $M=-I_2$.}

In this case, we have $(S_M^+(-1), S_M^-(-1))=(1, 1)$.
Thus as in Case 4, we have
\bea i_{-1}(\gamma^{2m-1})=(2m-1)i_1(\gamma)-1,\quad
\nu_{-1}(\gamma^{2m-1})=2.\lb{3.33}\eea

{\bf Case 6.} {\it $M$ is conjugate to a matrix
$\left(\matrix{ -1 & b\cr
                  0 & -1\cr}\right)$ for some $b>0$.}

In this case, we have $(S_M^+(-1), S_M^-(-1))=(0, 0)$.
Thus by Theorem 9.2.1 of \cite{Lon4}, we have
\bea i_{-1}(\gamma^{2m-1})&=&\sum_{\omega^{2m-1}=-1}i_\omega(\gamma)
=\sum_{k=1}^{2m-1}i_{\frac{(2k-1)\pi}{2m-1}}(\gamma)
=(2m-1)i_1(\gamma),\nn\\
\nu_{-1}(\gamma^{2m-1})&=&1.\lb{3.34}\eea

{\bf Case 7.} {\it
$M=\left(\matrix{\cos\th &-\sin\th\cr
                      \sin\th & \cos\th\cr}\right)$ with some
$\th\in (0,\pi)\cup (\pi,2\pi)$. }

In this case, we have $(S_M^+(e^{\sqrt{-1\theta}}), S_M^-(e^{\sqrt{-1\theta}}))=(0, 1)$.
Thus by Theorem 9.2.1 of \cite{Lon4} and Lemma 3.8, we have
\bea i_{-1}(\gamma^{2m-1})&=&\sum_{\omega^{2m-1}=-1}i_\omega(\gamma)
=\sum_{k=1}^{2m-1}i_{\frac{(2k-1)\pi}{2m-1}}(\gamma)
\nn\\&=&\sum_{2k-1<\frac{(2m-1)\theta}{\pi}}i_1(\gamma)
+\sum_{\frac{(2m-1)\theta}{\pi}\le2k-1\le\frac{(2m-1)(2\pi-\theta)}{\pi}}(i_1(\gamma)-1)
\nn\\&&+\sum_{\frac{(2m-1)(2\pi-\theta)}{\pi}<2k-1\le 4m-2}i_1(\gamma)
\nn\\&=&(2m-1)(i_1(\gamma)-1)+2E\left(\frac{(2m-1)\theta}{2\pi}+\frac{1}{2}\right)-2,
\nn\\
\nu_{-1}(\gamma^{2m-1})&=&2-2\phi\left(\frac{(2m-1)\theta}{2\pi}+\frac{1}{2}\right),\lb{3.35}\eea
provided
$\theta\in(0, \pi)$. When
$\theta\in(\pi, 2\pi)$, we have
\bea i_{-1}(\gamma^{2m-1})&=&\sum_{\omega^{2m-1}=-1}i_\omega(\gamma)
=\sum_{k=1}^{2m-1}i_{\frac{(2k-1)\pi}{2m-1}}(\gamma)
\nn\\&=&\sum_{2k-1\le\frac{(2m-1)(2\pi-\theta)}{\pi}}i_1(\gamma)
+\sum_{\frac{(2m-1)(2\pi-\theta)}{\pi}<2k-1<\frac{(2m-1)\theta}{\pi}}(i_1(\gamma)+1)
\nn\\&&+\sum_{\frac{(2m-1)\theta}{\pi}\le2k-1\le 4m-2}i_1(\gamma)
\nn\\&=&(2m-1)(i_1(\gamma)-1)+2E\left(\frac{(2m-1)\theta}{2\pi}+\frac{1}{2}\right)-2,
\nn\\
\nu_{-1}(\gamma^{2m-1})&=&2-2\phi\left(\frac{(2m-1)\theta}{2\pi}+\frac{1}{2}\right).\nn\eea

{\bf Case 8.} {\it $M=\left(\matrix{R(\th) & b
                  \cr 0 & R(\th)\cr}\right)$    with some
$\th\in (0,\pi)\cup (\pi,2\pi)$ and $b=
\left(\matrix{ b_1 &b_2\cr b_3 & b_4\cr}\right)\in\R^{2\times2}$,
such that $(b_2-b_3)\sin\theta<0$.}

In this case, we have $(S_M^+(e^{\sqrt{-1\theta}}), S_M^-(e^{\sqrt{-1\theta}}))=(1, 1)$.
Thus by Theorem 9.2.1 of \cite{Lon4}, we have
\bea i_{-1}(\gamma^{2m-1})&=&\sum_{\omega^{2m-1}=-1}i_\omega(\gamma)
=\sum_{k=1}^{2m-1}i_{\frac{(2k-1)\pi}{2m-1}}(\gamma)
\nn\\&=&(2m-1)i_1(\gamma)+2\phi\left(\frac{(2m-1)\theta}{2\pi}+\frac{1}{2}\right)-2,
\nn\\
\nu_{-1}(\gamma^{2m-1})&=&2-2\phi\left(\frac{(2m-1)\theta}{2\pi}+\frac{1}{2}\right).\lb{3.36}\eea

{\bf Case 9.} {\it $M=\left(\matrix{R(\th) & b
                  \cr 0 & R(\th)\cr}\right)$    with some
$\th\in (0,\pi)\cup (\pi,2\pi)$ and $b=
\left(\matrix{ b_1 &b_2\cr b_3 & b_4\cr}\right)\in\R^{2\times2}$,
such that $(b_2-b_3)\sin\theta>0$.}

In this case, we have $(S_M^+(e^{\sqrt{-1\theta}}), S_M^-(e^{\sqrt{-1\theta}}))=(0, 0)$.
Thus by Theorem 9.2.1 of \cite{Lon4}, we have
\bea i_{-1}(\gamma^{2m-1})&=&\sum_{\omega^{2m-1}=-1}i_\omega(\gamma)
=\sum_{k=1}^{2m-1}i_{\frac{(2k-1)\pi}{2m-1}}(\gamma)
=(2m-1)i_1(\gamma),
\nn\\
\nu_{-1}(\gamma^{2m-1})&=&2-2\phi\left(\frac{(2m-1)\theta}{2\pi}+\frac{1}{2}\right).\lb{3.37}\eea

{\bf Case 10.} {\it $M$ is hyperbolic, i.e., $\sigma(M)\cap\U=\emptyset$.}

In this case,  by Theorem 9.2.1 of \cite{Lon4}, we have
\bea i_{-1}(\gamma^{2m-1})&=&\sum_{\omega^{2m-1}=-1}i_\omega(\gamma)
=\sum_{k=1}^{2m-1}i_{\frac{(2k-1)\pi}{2m-1}}(\gamma)
=(2m-1)i_1(\gamma),\nn\\
\nu_{-1}(\gamma^{2m-1})&=&0.\lb{3.38}\eea

{\bf Proposition 3.11.} {\it
For any $m\in\N$, we have the estimate
\be i_{-1}(\gamma^{2m+1})-i_{-1}(\gamma^{2m-1})-\nu_{-1}(\gamma^{2m-1})
\ge 2i_1(\gamma)-e(M).
\lb{3.39}\ee}

{\bf Proof.} We consider each of the above cases.

{\bf Case 1.} {\it $M$ is conjugate to a matrix
$\left(\matrix{ 1 & b\cr
                  0 & 1\cr}\right)$ for some $b>0$.}

In this case we have
\bea i_{-1}(\gamma^{2m+1})-i_{-1}(\gamma^{2m-1})-\nu_{-1}(\gamma^{2m-1})
=2i_1(\gamma)+2.
\nn\eea

{\bf Case 2.} {\it $M=I_2$, the $2\times 2$ identity matrix}.

In this case we have
\bea i_{-1}(\gamma^{2m+1})-i_{-1}(\gamma^{2m-1})-\nu_{-1}(\gamma^{2m-1})
=2i_1(\gamma)+2.
\nn\eea

{\bf Case 3.} {\it $M$ is conjugate to a matrix
$\left(\matrix{ 1 & b\cr
                  0 & 1\cr}\right)$ for some $b<0$.}

In this case we have
\bea i_{-1}(\gamma^{2m+1})-i_{-1}(\gamma^{2m-1})-\nu_{-1}(\gamma^{2m-1})
=2i_1(\gamma).
\nn\eea

{\bf Case 4.} {\it $M$ is conjugate to a matrix
$\left(\matrix{ -1 & b\cr
                  0 & -1\cr}\right)$ for some $b<0$.}

In this case we have
\bea i_{-1}(\gamma^{2m+1})-i_{-1}(\gamma^{2m-1})-\nu_{-1}(\gamma^{2m-1})
=2i_1(\gamma)-1
\nn\eea

{\bf Case 5.} {\it $M=-I_2$, the $2\times 2$ identity matrix}.

In this case we have
\bea i_{-1}(\gamma^{2m+1})-i_{-1}(\gamma^{2m-1})-\nu_{-1}(\gamma^{2m-1})
=2i_1(\gamma)-2.
\nn\eea

{\bf Case 6.} {\it $M$ is conjugate to a matrix
$\left(\matrix{ -1 & b\cr
                  0 & -1\cr}\right)$ for some $b>0$.}

In this case we have
\bea i_{-1}(\gamma^{2m+1})-i_{-1}(\gamma^{2m-1})-\nu_{-1}(\gamma^{2m-1})
=2i_1(\gamma)-1.
\nn\eea

{\bf Case 7.} {\it
$M=\left(\matrix{\cos\th &-\sin\th\cr
                      \sin\th & \cos\th\cr}\right)$ with some
$\th\in (0,\pi)\cup (\pi,2\pi)$. }

In this case we have
\bea &&i_{-1}(\gamma^{2m+1})-i_{-1}(\gamma^{2m-1})-\nu_{-1}(\gamma^{2m-1})
\nn\\=&&2(i_1(\gamma)-1)+2E\left(\frac{(2m+1)\theta}{2\pi}+\frac{1}{2}\right)
-2E\left(\frac{(2m-1)\theta}{2\pi}+\frac{1}{2}\right)
\nn\\&&-\left(2-2\phi\left(\frac{(2m-1)\theta}{2\pi}+\frac{1}{2}\right)\right)
\nn\\\ge&&2(i_1(\gamma)-1).
\nn\eea

{\bf Case 8.} {\it $M=\left(\matrix{R(\th) & b
                  \cr 0 & R(\th)\cr}\right)$    with some
$\th\in (0,\pi)\cup (\pi,2\pi)$ and $b=
\left(\matrix{ b_1 &b_2\cr b_3 & b_4\cr}\right)\in\R^{2\times2}$,
such that $(b_2-b_3)\sin\theta<0$.}

In this case we have
\bea &&i_{-1}(\gamma^{2m+1})-i_{-1}(\gamma^{2m-1})-\nu_{-1}(\gamma^{2m-1})
\nn\\=&&2i_1(\gamma)+2\phi\left(\frac{(2m+1)\theta}{2\pi}+\frac{1}{2}\right)
-2\phi\left(\frac{(2m-1)\theta}{2\pi}+\frac{1}{2}\right)
\nn\\&&-\left(2-2\phi\left(\frac{(2m-1)\theta}{2\pi}+\frac{1}{2}\right)\right)
\nn\\\ge&&2i_1(\gamma)-2.
\nn\eea

{\bf Case 9.} {\it $M=\left(\matrix{R(\th) & b
                  \cr 0 & R(\th)\cr}\right)$    with some
$\th\in (0,\pi)\cup (\pi,2\pi)$ and $b=
\left(\matrix{ b_1 &b_2\cr b_3 & b_4\cr}\right)\in\R^{2\times2}$,
such that $(b_2-b_3)\sin\theta>0$.}

In this case we have
\bea &&i_{-1}(\gamma^{2m+1})-i_{-1}(\gamma^{2m-1})-\nu_{-1}(\gamma^{2m-1})
\nn\\=&&2i_1(\gamma)-\left(2-2\phi\left(\frac{(2m-1)\theta}{2\pi}+\frac{1}{2}\right)\right)
\nn\\\ge&&2i_1(\gamma)-2.
\nn\eea

{\bf Case 10.} {\it $M$ is hyperbolic, i.e., $\sigma(M)\cap\U=\emptyset$.}

In this case we have
\bea i_{-1}(\gamma^{2m+1})-i_{-1}(\gamma^{2m-1})-\nu_{-1}(\gamma^{2m-1})
=2i_1(\gamma).
\nn\eea
Combining the above cases, we obtain the proposition.\hfill\hb

\setcounter{equation}{0}
\section{Proof of the main theorem}

In this section, we give the proof of the main theorem.
first we have the following.

{\bf Lemma 4.1.} {\it Suppose $u^{2k-1}$ is a nonzero critical
point of $\Psi_a$ such that $u$ corresponds to $(\tau, y)\in \mathcal{T}_s(\Sg)$.
Then we can find $m\in\N$ such that
\be i(u^{2m+1})-i(u^{2m-1})\ge 4.
\lb{4.1}\ee}

{\bf Proof.} Let $(\tau,\, y)\in\mathcal{T}_s(\Sigma)$. The fundamental
solution $\gamma_y : [0,\,\tau/2]\rightarrow \Sp(2n)$ with $\gamma_y(0)=I_{2n}$
of the linearized Hamiltonian system
\be \dot w(t)=JH^{\prime\prime}(y(t))w(t),\qquad \forall t\in\R,\lb{4.2}\ee
is called the {\it associate symplectic path} of $(\tau,\, y)$.
Then as in \S1.7 of \cite{Eke3}, we have
\bea \gamma_y(\tau/2)=\left(\matrix{-I_2\quad 0  \cr
                        0\quad\;\; C   \cr }\right)\lb{4.3}
\eea
in an appropriate coordinate. Then by Lemma 3.1 and Theorem 3.5,
we have
\be i(u^{2k-1})= i_{-1}(\gamma^{2k-1}),\quad \nu(u^{2k-1})=\nu_{-1}(\gamma^{2k-1}),
      \lb{4.4}\ee
for any $k\in\N$. By Theorem 3.10, the matrix $\gamma_y(\tau/2)$ can be connected
in $\Om^0(\gamma_y(\tau/2))$ to a basic form decomposition
$M= (-I_2)\dm M_1\dm\cdots\dm M_l$.
Since $n\ge 2$, we may write $M=(-I_2)\dm M_1\dm M^\prime$,
where $M^\prime=M_2\dm\cdots\dm M_l$.
By the symplectic additivity of indices, cf. \cite{Lon2}-\cite{Lon4},
 we have
\be i_{-1}(\gamma^{2k-1})=i_{-1}(\gamma_1^{2k-1})+i_{-1}(\gamma_2^{2k-1})
\lb{4.3}\ee
where $\gamma_1$ and $\gamma_2$ are appropriate symplectic
paths such that $\gamma_1(\tau/2)=(-I_2)\dm M_1$ and
$\gamma_2(\tau/2)=M^\prime$.

Note that by Theorem 3.5, we have $i_1(\gamma)\ge n$.
Now we consider each case as in \S3.

{\bf Case 1.} {\it $M_1=\left(\matrix{ 1 & b\cr
                  0 & 1\cr}\right)$ for some $b>0$
or $M_1=I_2$.}

In this case we have
\bea  &&i(u^{2m+1})-i(u^{2m-1})
=i_{-1}(\gamma^{2m+1})-i_{-1}(\gamma^{2m-1})
\nn\\=&&i_{-1}(\gamma_1^{2m+1})-i_{-1}(\gamma_1^{2m-1})
+i_{-1}(\gamma_2^{2m+1})-i_{-1}(\gamma_2^{2m-1})
\nn\\\ge&&2i_1(\gamma_1)+2+2i_1(\gamma_2)-(2n-4)+\nu_{-1}(\gamma_2^{2m-1})
\nn\\\ge&&2i_1(\gamma)+6-2n\ge 6.
\nn\eea
Note that in the above computations, we use
(\ref{3.29}), (\ref{3.30}), (\ref{3.33}),
Proposition 3.11 and $i_1(\gamma)\ge n$.

{\bf Case 2.} {\it $M$ is conjugate to a matrix
$\left(\matrix{ 1 & b\cr
                  0 & 1\cr}\right)$ for some $b<0$.}

In this case, by (\ref{3.31}) we have
\bea  &&i(u^{2m+1})-i(u^{2m-1})
=i_{-1}(\gamma^{2m+1})-i_{-1}(\gamma^{2m-1})
\nn\\=&&i_{-1}(\gamma_1^{2m+1})-i_{-1}(\gamma_1^{2m-1})
+i_{-1}(\gamma_2^{2m+1})-i_{-1}(\gamma_2^{2m-1})
\nn\\\ge&&2i_1(\gamma_1)+2i_1(\gamma_2)-(2n-4)+\nu_{-1}(\gamma_2^{2m-1})
\nn\\\ge&&2i_1(\gamma)+4-2n\ge 4.
\nn\eea

{\bf Case 3.} {\it $M=\left(\matrix{ -1 & b\cr
                  0 & -1\cr}\right)$ for $b\in\R$.}

In this case, by (\ref{3.32})-(\ref{3.34}) we have
\bea  &&i(u^{2m+1})-i(u^{2m-1})
=i_{-1}(\gamma^{2m+1})-i_{-1}(\gamma^{2m-1})
\nn\\=&&i_{-1}(\gamma_1^{2m+1})-i_{-1}(\gamma_1^{2m-1})
+i_{-1}(\gamma_2^{2m+1})-i_{-1}(\gamma_2^{2m-1})
\nn\\\ge&&2i_1(\gamma_1)+2i_1(\gamma_2)-(2n-4)+\nu_{-1}(\gamma_2^{2m-1})
\nn\\\ge&&2i_1(\gamma)+4-2n\ge 4.
\nn\eea

{\bf Case 4.} {\it
$M=\left(\matrix{\cos\th &-\sin\th\cr
                      \sin\th & \cos\th\cr}\right)$ with some
$\th\in (0,\pi)\cup (\pi,2\pi)$. }

In this case, by (\ref{3.35}) we have
\bea  &&i(u^{2m+1})-i(u^{2m-1})
=i_{-1}(\gamma^{2m+1})-i_{-1}(\gamma^{2m-1})
\nn\\=&&i_{-1}(\gamma_1^{2m+1})-i_{-1}(\gamma_1^{2m-1})
+i_{-1}(\gamma_2^{2m+1})-i_{-1}(\gamma_2^{2m-1})
\nn\\\ge&&2i_1(\gamma_1)-2+2E\left(\frac{(2m+1)\theta}{2\pi}+\frac{1}{2}\right)
-2E\left(\frac{(2m-1)\theta}{2\pi}+\frac{1}{2}\right)
\nn\\&&+2i_1(\gamma_2)-(2n-4)+\nu_{-1}(\gamma_2^{2m-1})
\nn\\\ge&&2i_1(\gamma)+4-2n\ge 4
\nn\eea
provided we choose $m$ such that
$E\left(\frac{(2m+1)\theta}{2\pi}+\frac{1}{2}\right)
-E\left(\frac{(2m-1)\theta}{2\pi}+\frac{1}{2}\right)\ge 1$.

{\bf Case 5.} {\it $M=\left(\matrix{R(\th) & b
                  \cr 0 & R(\th)\cr}\right)$    with some
$\th\in (0,\pi)\cup (\pi,2\pi)$ and $b=
\left(\matrix{ b_1 &b_2\cr b_3 & b_4\cr}\right)\in\R^{2\times2}$,
such that $(b_2-b_3)\sin\theta<0$.}

In this case, by (\ref{3.36}) we have
\bea  &&i(u^{2m+1})-i(u^{2m-1})
=i_{-1}(\gamma^{2m+1})-i_{-1}(\gamma^{2m-1})
\nn\\=&&i_{-1}(\gamma_1^{2m+1})-i_{-1}(\gamma_1^{2m-1})
+i_{-1}(\gamma_2^{2m+1})-i_{-1}(\gamma_2^{2m-1})
\nn\\\ge&&2i_1(\gamma_1)+2\varphi\left(\frac{(2m+1)\theta}{2\pi}+\frac{1}{2}\right)
-2\varphi\left(\frac{(2m-1)\theta}{2\pi}+\frac{1}{2}\right)
\nn\\&&+2i_1(\gamma_2)-(2n-6)+\nu_{-1}(\gamma_2^{2m-1})
\nn\\\ge&&2i_1(\gamma)+4-2n\ge 4.
\nn\eea

{\bf Case 6.} {\it $M=\left(\matrix{R(\th) & b
                  \cr 0 & R(\th)\cr}\right)$    with some
$\th\in (0,\pi)\cup (\pi,2\pi)$ and $b=
\left(\matrix{ b_1 &b_2\cr b_3 & b_4\cr}\right)\in\R^{2\times2}$,
such that $(b_2-b_3)\sin\theta>0$.}

In this case, by (\ref{3.37}) we have
\bea  &&i(u^{2m+1})-i(u^{2m-1})
=i_{-1}(\gamma^{2m+1})-i_{-1}(\gamma^{2m-1})
\nn\\=&&i_{-1}(\gamma_1^{2m+1})-i_{-1}(\gamma_1^{2m-1})
+i_{-1}(\gamma_2^{2m+1})-i_{-1}(\gamma_2^{2m-1})
\nn\\\ge&&2i_1(\gamma_1)+2i_1(\gamma_2)-(2n-6)+\nu_{-1}(\gamma_2^{2m-1})
\nn\\\ge&&2i_1(\gamma)+6-2n\ge 6.
\nn\eea

{\bf Case 7.} {\it $M$ is hyperbolic, i.e., $\sigma(M)\cap\U=\emptyset$.}

In this case, by (\ref{3.38}) we have
\bea  &&i(u^{2m+1})-i(u^{2m-1})
=i_{-1}(\gamma^{2m+1})-i_{-1}(\gamma^{2m-1})
\nn\\=&&i_{-1}(\gamma_1^{2m+1})-i_{-1}(\gamma_1^{2m-1})
+i_{-1}(\gamma_2^{2m+1})-i_{-1}(\gamma_2^{2m-1})
\nn\\\ge&&2i_1(\gamma_1)+2i_1(\gamma_2)-(2n-4)+\nu_{-1}(\gamma_2^{2m-1})
\nn\\\ge&&2i_1(\gamma)+4-2n\ge 4.
\nn\eea
Combining all the above cases, we obtain the lemma.\hfill\hb

{\bf Proof of Theorem 1.1.}
We prove by contraction. Assume $\mathcal{T}_s(\Sg)=\{(\tau, y)\}$.
 Suppose $u^{2m-1}$ is a nonzero critical
point of $\Psi_a$ such that $u$ corresponds to $(\tau, y)\in \mathcal{T}_s(\Sg)$.
By Lemma 4.1, we may assume
$i(u^{2m+1})-i(u^{2m-1})\ge 4$.
The index interval of $(\tau, y)$ at $2m-1$ is defined to be
$\mathcal{G}_{2m-1}=(i(u^{2m-3})+\nu(u^{2m-3})-1,\;  i(u^{2m+1}))$.
Note that by Proposition 3.11 and $i_1(\gamma)\ge n$, we have
$i(u^{2m-3})+\nu(u^{2m-3})\le i(u^{2m-1})$.
Thus we have $(i(u^{2m-1})-1,\;  i(u^{2m+1}))\subset\mathcal{G}_{2m-1}$.
Hence we can find two distinct even integers $2T_1, 2T_2\in\mathcal{G}_{2m-1}$.
Let $c_{T_1+1}$ and $c_{T_2+1}$ be the two critical values of
$\Psi_a$ found by Proposition 2.16. Then we have
$c_{T_1+1}\neq c_{T_2+1}$ since $^\#\mathcal{T}_s(\Sg)<\infty$.
By Proposition 2.17, we have
\bea \Psi_a(u^{2m_1-1})=c_{T_1+1},\qquad i(u^{2m_1})\le 2T_1\le i(u^{2m_1-1}) + \nu(u^{2m_1-1})-1,
\nn\\\Psi_a(u^{2m_2-1})=c_{T_2+1},\qquad i(u^{2m_2})\le 2T_2\le i(u^{2m_2-1}) + \nu(u^{2m_2-1})-1,
\eea
for some $m_1, m_2\in\N$.
On the other hand, we must have $m_1=m_2$ by Proposition 3.11.
Thus we have $c_{T_1+1}=c_{T_2+1}$. This contradiction proves the theorem.
\hfill\hb

\noindent {\bf Acknowledgements.} I would like to sincerely thank my
Ph. D. thesis advisor, Professor Yiming Long, for introducing me to Hamiltonian
dynamics and for his valuable help and encouragement during my research.
I would like to say that how enjoyable it is to work with him.

\bibliographystyle{abbrv}

\begin{thebibliography}{99}

\bibitem[Bor1]{Bor1} A. Borel, Seminar on Transformation Groups. Princeton
Univ. Press. Princeton. 1960.
\bibitem[CoZ1]{CoZ1} C. Conley and E. Zehnder, Morse-type index theory for
flows and periodic solutions for Hamiltonian equations.
{\it Comm. Pure. Appl. Math.} 37 (1984) 207-253.
\bibitem[Eke1]{Eke1} I. Ekeland, Une th\'{e}orie de Morse pour les
syst\`{e}mes hamiltoniens convexes. {\it Ann. IHP. Anal. non
Lin\'{e}aire.} 1 (1984) 19-78.
\bibitem[Eke2]{Eke2} I. Ekeland, An index throry for periodic solutions
of convex Hamiltonian systems. {\it Proc. Symp. in Pure Math.} 45 (1986)
395-423.
\bibitem[Eke3]{Eke3} I. Ekeland, Convexity Methods in Hamiltonian
Mechanics. Springer-Verlag. Berlin. 1990.
\bibitem[EkH1]{EkH1} I. Ekeland and H. Hofer,  Convex Hamiltonian energy
surfaces and their closed trajectories.  {\it Comm. Math. Phys.}
113 (1987) 419-467.
\bibitem[EkL1]{EkL1} I. Ekeland and L. Lassoued,  Multiplicit\'e des
trajectoires ferm\'ees d'un syst\'eme hamiltonien sur une
hypersurface d'energie convexe. {\it Ann. IHP. Anal. non Lin\'eaire}.
4 (1987) 1-29.
\bibitem[FaR1]{FaR1} E. Fadell and P. Rabinowitz, Generalized
comological index throries for Lie group actions with an application
to bifurcation equations for Hamiltonian systems.
{\it Invent. Math. } 45 (1978) 139-174.
\bibitem[GrM1]{GrM1} D. Gromoll and W. Meyer,  On differentiable functions
with isolated critical points. {\it Topology.} 8 (1969) 361-369.
\bibitem[GrM2]{GrM2} D. Gromoll and W. Meyer,  Periodic geodesics on compact
Riemannian manifolds. {\it J. Diff. Geod.} 3 (1969) 493-510.
\bibitem[HWZ]{HWZ} H. Hofer, K. Wysocki, and E. Zehnder, The dynamics on
three-dimensional strictly convex energy surfaces. {\it Ann. of
Math.} 148 (1998) 197-289.
\bibitem[Lon1]{Lon1} Y. Long, Maslov-type index, degenerate critical points
and asymptotically linear Hamiltonian systems.
{\it Science in China.} Series A. 33(1990), 1409-1419.
\bibitem[Lon2]{Lon2} Y. Long,  Bott formula of the Maslov-type index
theory. {\it Pacific J. Math.} 187 (1999), 113-149.
\bibitem[Lon3]{Lon3} Y. Long,  Precise iteration formulae of the
Maslov-type index theory and ellipticity of closed characteristics.
{\it Advances in Math.} 154 (2000), 76-131.
\bibitem[Lon4]{Lon4} Y. Long,  Index Theory for Symplectic Paths with
Applications. Progress in Math. 207, Birkh\"auser. Basel. 2002.
\bibitem[Lon5]{Lon5} Y. Long, Index iteration theory for symplectic paths
with applications to nonlinear Hamiltonian systems. {\it Proc. of Inter.
Congress of Math. 2002.} Vol.II, 303-313. Higher Edu. Press. Beijing. 2002.
\bibitem[Lon6]{Lon6} Y. Long,  Index iteration theory for symplectic paths and
multiple periodic solution orbits. {\it Frontiers of Math. in China.} 1 (2006)
178-201.
\bibitem[LLZ]{LLZ} C. Liu, Y. Long and C. Zhu, Multiplicity of closed
characteristics on symmetric convex hypersurfaces in $\R^{2n}$.
{\it Math. Ann.} 323 (2002), 201-215.
\bibitem[LZe1]{LZe1} Y. Long and E. Zehnder, Morse theory for forced
oscillations of asymptotically linear Hamiltonian systems.
In {\it Stoc. Proc. Phys. and Geom.}, S. Albeverio et al. ed.
World Sci. (1990) 528-563.
\bibitem[LoZ1]{LoZ1} Y. Long and C. Zhu,  Closed characteristics on
compact convex hypersurfaces in $\R^{2n}$.  {\it Ann. of Math.}
155 (2002) 317-368.
\bibitem[Rab1]{Rab1} P. H. Rabinowitz,  Periodic solutions of Hamiltonian
systems. {\it Comm. Pure Appl. Math}. 31 (1978) 157-184.
\bibitem[Szu1]{Szu1} A. Szulkin, Morse theory and existence of periodic solutions
of convex Hamiltonian systems. {\it Bull. Soc. Math. France.} 116 (1988) 171-197.
\bibitem[Vit1]{Vit1} C. Viterbo,   Equivariant Morse theory
for starshaped Hamiltonian systems. {\it Trans. Amer. Math. Soc. }
311 (1989) 621-655.
\bibitem[Vit2]{Vit2} C. Viterbo,  A new obstruction to embedding Lagrangian tori.
{\it Invent. Math.} 100 (1990) 301-320.
\bibitem[W1]{W1} W. Wang, Stability of closed characteristics on compact convex
hypersurfaces in $R^6$. {\it J. Eur. Math. Soc.} 11 (2009), 575-596.
\bibitem[W2]{W2} W. Wang, Symmetric closed characteristics on symmetriccompact
convex hypersurfaces in $R^{2n}$,
{\it J.  Diff. Equa.} 246 (2009), 4322-4331.
\bibitem[WHL]{WHL} W. Wang, X. Hu and Y. Long, Resonance identity, stability and
multiplicity of closed characteristics on compact convex
hypersurfaces.  {\it  Duke Math. J.} 139 (2007) 411-462.
\bibitem[Was1]{Was1} A. Wasserman,  Morse theory for G-manifolds.
 {\it Bull. Amer. Math. Soc. } (March, 1965).
\bibitem[Wei1]{Wei1} A. Weinstein, Periodic orbits for convex Hamiltonian
systems. {\it Ann. of Math}. 108 (1978) 507-518.

\end{thebibliography}

\medskip

\end{document}